%% file: main.tex
\documentclass[12pt,a4paper,preprint,sort&compress]{elsarticle}

\usepackage[dvipsnames]{xcolor}
\usepackage{graphicx}
\usepackage{amssymb}
\usepackage{amsmath,amsthm}
\usepackage{bm}
\usepackage{subfigure}
\usepackage[T1]{fontenc}
\usepackage{amsthm}
\usepackage{tikz}
\usepackage{pgfplots}
\usepackage{etoolbox}
\usepgfplotslibrary{patchplots}
\usepgfplotslibrary{colormaps}
\usepgfplotslibrary{groupplots}
\usepackage{algorithm2e}
\usepackage{soul}
\usepackage{stmaryrd}
\usepackage{enumerate}

\theoremstyle{plain}
\newtheorem{assumption}{Assumption}
\newtheorem{remark}{Remark}

\makeatletter
\def\ps@pprintTitle{%
 \let\@oddhead\@empty
 \let\@evenhead\@empty
 \def\@oddfoot{}%
 \let\@evenfoot\@oddfoot}
\makeatother

\usepackage[hidelinks]{hyperref}

\input{macros}

\newif\ifrecompiletikz
\recompiletikzfalse

\ifrecompiletikz
\usetikzlibrary{positioning,arrows,mindmap,backgrounds, shapes, decorations.fractals, calc, intersections,decorations.pathmorphing,patterns}
 \usetikzlibrary{external}
 \tikzexternaldisable
\fi

	\title{Fast assembly of Galerkin matrices for 3D solid laminated composites using finite element and isogeometric discretizations}
	\author[lausanne]{Pablo Antolin\corref{cor1}}

	\cortext[cor1]{Corresponding author,  E-mail: \texttt{pablo.antolin@epfl.ch}}
	\address[lausanne]{Institute of Mathematics - \'Ecole Polytechnique F\'ed\'erale de Lausanne\\
		CH-1015 Lausanne, Switzerland}
	\date{\empty}

\begin{document}
		
	\begin{abstract}
	  This work presents a novel methodology for speeding up the assembly of stiffness matrices for laminate composite 3D structures in the context of isogeometric and finite element discretizations.
		By splitting the involved terms into their in-plane and out-of-plane contributions, this method computes the problems's 3D stiffness matrix as a combination of 2D (in-plane) and 1D (out-of-plane) integrals.
		Therefore, the assembly's computational complexity is reduced to the one of a 2D problem.
		Additionally, the number of 2D integrals to be computed becomes independent of the number of material layers that constitute the laminated composite, it only depends on the number of different materials used (or different orientations of the same anisotropic material). 
		Hence, when a high number of layers is present, the proposed technique reduces by orders of magnitude the computational time required to create the stiffness matrix with standard methods, being the resulting matrices identical up to machine precision.
		The predicted performance is illustrated through numerical experiments.
	\end{abstract}
	\begin{keyword}
  Composite laminates; Finite elements; Isogeometric analysis; Fast matrix assembly; Layerwise theory
	\end{keyword}

	\maketitle
  \setcounter{tocdepth}{1}


\section{Introduction} \label{sec:intro}
Laminated composite materials are widely used in the aerospace and automotive as well as construction industries.
Their lightness, strength and stiffness make them a very appealing solution in multiple situations, and those are also the reasons behind their continuous spread to many other consumer goods during the last decades.

Composite laminates are made of stacked material layers, also denoted as plies, that present different mechanical properties.
These laminae are frequently made of stiff reinforcement fibers embedded in a soft matrix material that glues them together, as it is the case of glass or carbon fiber materials.
The constituent plies are usually stacked in such a way that their main fiber direction follows different orientations.
Classical stack-up sequences are the $0^\circ/90^\circ/0^\circ/90^\circ/\dots$ and $0^\circ/\pm45^\circ/90^\circ/\dots$ layer setups, in which plies of the same material are stacked using 2 or 4 different in-plane fiber orientation angles, respectively (see, e.g., \cite{reddy_practical_1995}).

Due to this heterogeneous stack of layers, laminated structures present a quite complex mechanical response when loaded.
In order to correctly assess the possible failure modes it is quite important to characterize precisely the strain and stress distributions in the different material layers and at their interfaces.
We refer the interested reader to the classical references \cite{gibson_principles_2016,jones_mechanics_2018,vinson_behavior_2004,reddy_mechanics_2004} for a deeper insight on the many mechanical aspects of composite materials.

Beyond academic cases for which closed form solutions exist (e.g., \cite{pagano_exact_1970,varadan_bending_1991}), the complexity of the laminates elastic behavior calls for the use of numerical methods.
And, without doubt, the finite element method (FEM) is nowadays, and has been during the last decades, the most popular tool for the analysis of their mechanical response.
A review of the different finite element methods available for their study is out of the scope of this paper, but we refer the interested reader to the numerous reviews on the topic, e.g. \cite{carrera_theories_2002,carrera_theories_2003,reddy_theories_1994,carrera_survey_2008,zhang_recent_2009,s_qatu_review_2012,kreja_literature_2011,liew_overview_2019}, and the many references therein.

From that vast literature it is easy to realize that a high percentage of the published finite element models correspond to 2D plate and shell models \cite{zhang_recent_2009,kreja_literature_2011,liew_overview_2019}.
The two-dimensional FEM models are primarily founded on two main composite theories: the equivalent single-layer \cite{kreja_literature_2011,kulkarni_review_2018} and layerwise \cite{reddy_mechanics_2004,carrera_multilayered_1999,ferreira_analysis_2005} approaches.
As discussed in the review \cite{liew_overview_2019}, layerwise based FEM models, that present good compromise between accuracy and computational cost, combine an in-plane 2D FEM discretization with additional degrees of freedom for describing the displacements and/or transverse stresses of all the material layers.
In these approaches the number of required degrees of freedom scales linearly with the number of plies.

Similarly, 3D elasticity models are able to accurately reproduce the laminate's strain and stress profiles by representing each ply with one or more elements along its thickness.
As in the case of layerwise approaches, the number of degrees of freedom scales linearly with the number of layers, but, additionally, the computational cost of the matrix assembly is very high compared to 2D models \cite{liew_overview_2019}.
In this paper we propose a new methodology that significantly reduces this assembly cost, becoming analogous to the complexity of a 2D problem and independent of the total number of material layers.

The introduction of the Isogeometric Analysis (IGA) concept by Hughes et al. in \cite{hughes_isogeometric_2005,cottrell_isogeometric_2009} constituted a great success.
Based on the high continuity of spline basis functions and their superior approximation properties compared to classical finite elements methods \cite{evans_n-widths_2009,beirao_da_veiga_estimates_2011}, IGA has proven to be a powerful tool for a broad spectrum of applications:
from solid mechanics \cite{cottrell_isogeometric_2006,elguedj_$bar_2008,kiendl_isogeometric_2009} to fluid dynamics \cite{bazilevs_variational_2007,hsu_dynamic_2015} and other fields \cite{gomez_isogeometric_2008,buffa_isogeometric_2010}.
This success is also extensive to the analysis of laminated composites, where different authors proposed reduced dimension methods for plates, shells and beams \cite{bazilevs_3d_2010,kapoor_interlaminar_2013,nguyen-xuan_isogeometric_2013,thai_isogeometric_2015,remmers_isogeometric_2015,pavan_bending_2017,farzam_new_2018,ghafari_isogeometric_2019}.

On the other hand, despite their high accuracy, 3D solid IGA models in which at least one element along the thickness is used for every material ply \cite{guo_layerwise_2014,guo_layerwise_2015}, present a reduced interest due to their very high computational cost, as in the case of similar FEM approaches.
Even if in practice IGA models are able to achieve the same level of accuracy as FEM models, but using much coarser in-plane discretizations, the assembly cost per degree of freedom is often a major bottleneck.
To reduce this cost is a quite active research field \cite{hughes_efficient_2010,auricchio_simple_2012,schillinger_reduced_2014,antolin_efficient_2015,barton_optimal_2016,barton_gaussgalerkin_2017,calabro_fast_2017,mantzaflaris_low_2017,hiemstra_fast_2019}.

In order to overcome this burden, in \cite{dufour_cost-effective_2018} we proposed the use of 3D isogeometric models with a single element through the full thickness.
In a post-processing stage a very high-fidelity stress profile is recovered from the coarse simulation result.
This approach has been recently extended in \cite{patton_fast_2019} to the case of collocation methods.
In \cite{dufour_cost-effective_2018}, despite the use of single element through the thickness, the stiffness matrix assembly cost is still high, due to the use of layerwise quadrature schemes along the lamina thickness (the number of quadrature points scales linearly with the number of layers).

In this work, inspired by the sum-factorization technique \cite{melenk_fully_2001,ainsworth_bernsteinbezier_2011,antolin_efficient_2015}, we propose a new method for the assembly of 3D laminated composites for both finite element and isogeometric methods.
Our approach reduces drastically their assembly cost and is very effective for laminated structures that present many plies with the same configuration repeated along the laminate stack-up sequence (e.g., for $0^\circ/90^\circ/0^\circ/90^\circ/\dots$ and $0^\circ/\pm45^\circ/90^\circ/\dots$ layer setups).
By decomposing the computation of the required 3D integrals into its in-plane and out-of-plane contributions, the complexity of our assembly procedure scales as the assembly cost of a 2D problem, multiplied by the number of different ply configurations (for instance, 2 configurations for $0^\circ/90^\circ/0^\circ/90^\circ/\dots$ and 4 for $0^\circ/\pm45^\circ/90^\circ/\dots$), being independent of the total number of laminae.

The computed stiffness matrices with this new technique are identical, up to machine precision, to the ones obtained using a standard assembly procedure.

The assembly times obtained with the proposed method are even faster than the ones achieved with the standard method, using a simplified 3D model with a single element through the whole lamina's thickness, and combined with a material homogenization approach (e.g., \cite{sun_three-dimensional_1988}).
In addition, the homogenization approach would lead to less accurate results.

The rest of this paper is structured as follows.
Initially, Galerkin methods for the analysis of laminate composite structures are introduced in Section \ref{sec:galerkin}, using finite element and isogeometric discretizations.
The proposed fast assembly method is presented in Section \ref{sec:fast_assembly}, as well as a discussion of its theoretical computational complexity.
In Section \ref{sec:numerical} numerical experiments that consider different isogeometric discretizations and laminate configurations illustrate the performance of the proposed method confronted with classical assembly techniques.
Conclusions are drawn in Section \ref{sec:conclusions}. Finally, in \ref{ap:alternative} an alternative implementation that overcomes the use of Voigt's notation is detailed.
		
\section{Galerkin methods for 3D laminated composites} \label{sec:galerkin}

In this section we begin by introducing the linear elasticity problem for 3D composite laminates in the context of finite element and isogeometric discretizations.
We initially describe the continuous problem and its notation in Section \ref{sec:continuous}, followed by the family of used discretizations in Section \ref{sec:discrete}, that we particularize, in Section \ref{sec:multi-layered}, to the case of multi-layered materials.
Finally, the stiffness matrix's computational complexity is discussed in Section \ref{sec:classical_assembly}.

\subsection{Continuous elasticity problem} \label{sec:continuous}
Let us start by setting the continuous three-dimensional elasticity problem in strong form, namely
\begin{align} \label{eq:strong_problem}
\begin{array}{rll}
\nabla\cdot\stress &= \bm{f} &\text{in }\Omega\,,\\
\bm{u} &= \bm{0} &\text{on }\partial\Omega\,,
\end{array}
\end{align}
where $\bm{u}\in\mathbb{R}^3$ is displacement of the elastic body at every point, i.e., the problem unknown, $\bm{f}\in\mathbb{R}^3$ is the applied loading and $\strain\in\mathbb{R}^{3\times3}$ the small strain tensor, computed as  $\strain(\bm{u}) = \nabla^s\bm{u}$.
For small strains and linear materials the stress tensor $\stress\in\mathbb{R}^{3\times3}$ is computed as $\stress=\CC:\strain(\bm{u})$, where $\CC\in\mathbb{R}^{3\times3\times3\times3}$ is the fourth order elasticity tensor.
Without constituting any limitation, and for the sake of simplicity, only homogeneous Dirichlet boundary conditions are considered.

The domain $\Omega\subset\mathbb{R}^3$ occupied by the elastic body is the image of the parametric domain $\hat \Omega=[0,1]^3$ mapped with $\bm{F}:\hat\Omega\to\Omega$, which we assume to be a bi-Lipschitz homeomorphism.
Thus, a parametric point $\bm\xi=(\xi^1, \xi^2, \xi^3)$ is mapped into its physical image $\bm{x}$ as $\bm{F}:\bm{\xi}\in\hat\Omega\mapsto\bm{x}\in\Omega$.

In a classical way, considering the functional space $\bm{V} := \{\bm{u}\in H^1(\Omega)^3\, : \, \left.\bm{u}\right|_{\partial\Omega} = \bm{0}\}$,
the variational form associated to the strong problem \eqref{eq:strong_problem} can be formulated as: find $\bm{u}\in \bm{V}$ such that:
\begin{align} \label{eq:weak_problem}
a(\bm{v}, \bm{u}) = \int_\Omega \bm{v}\cdot\bm{f}\diff\bm{x},\quad\forall \bm{v}\in \bm{V}\,,
\end{align}
where the bilinear for $a:\bm{V}\times \bm{V}\to\mathbb{R}$ is:
\begin{align}\label{eq:a}
  a(\bm{v},\,\bm{u}) &= \int_\Omega \strain(\bm{v}):\mathbb{C}:\strain(\bm{u})\diff\bm{x}\,.
\end{align}

\subsection{Discrete elasticity problem}\label{sec:discrete}

In order to discretize the strong problem \eqref{eq:weak_problem}, we follow a standard Galerkin approach and choose a finite dimensional space $\bm{V}_h\subset \bm{V}$, such that $\bm{V}_h=(V_h)^3$ and:
\begin{subequations}
\begin{align}
&V_h = \{\hat u_h\circ \bm{F}^{-1} : \hat u_h \in N(\hat \Omega_h),  \, \left.\hat u_h\right|_{\partial\hat\Omega} = 0\}\,,\\
&N(\hat\Omega_h)  = \text{span}\left\lbrace \hat B^i(\bm \xi),\, i=1,\dots,n\right\rbrace\,,
\end{align}
\end{subequations}
where $\hat B^i(\bm \xi)$ are the space basis functions, $n$ is the space dimension and $\hat\Omega_h$ is a generic partition of the parametric domain.

\begin{assumption} \label{assmpt:0}
In this work, the basis functions $\hat B^i(\bm \xi)$ are chosen such that
\begin{align} \label{eq:basis_functions}
\hat B^i(\bm \xi)=\hat S^{i_s}(\xi^1,\xi^2)\,\hat T^{i_t}(\xi^3)\,,
\end{align}
with $i=(i_t - 1)\,n_s + i_s$, $i_s=1,\dots,n_s$ and $i_t=1,\dots,n_t$.
I.e., the basis functions $\hat B^i$ are created as the combination of $n_s$ in-plane functions $\hat S$ and $n_t$ out-of-plane functions $\hat T$.
Therefore, the space dimension is $n=n_s n_t$.
\end{assumption} 

Finite element spaces composed of Lagrangian hexahedron and wedge elements carry basis functions that fall into the category defined in \eqref{eq:basis_functions}.
However, this is not the case of tetrahedral finite element meshes.
We refer the interested reader to the classical references \cite{hughes_finite_1987,zienkiewicz_finite_2013,brenner_mathematical_2008,ciarlet_finite_2002} for further details on the definition of finite element spaces.

Regarding isogeometric discretizations, the basis functions definition \eqref{eq:basis_functions}
allows to consider different in-plane discretizations combined with non-rational B-spline basis functions along the third direction.
Thus, standard B-spline basis functions \cite{hughes_isogeometric_2005,cottrell_isogeometric_2009}, including NURBS, can be used in-plane, but also non-tensor product schemes such as HR-splines \cite{giannelli_thb-splines:_2012,vuong_hierarchical_2011}, LR-splines \cite{dokken_polynomial_2013,bressan_properties_2013} or T-splines \cite{beirao_da_veiga_analysis-suitable_2013,bazilevs_isogeometric_2010}.

Thus, based on the above defined finite space $\bm{V}_h$, we discretize the trial $\bm u$ and test  $\bm v$ functions of the problem \eqref{eq:weak_problem} by means of
\begin{align*}
\bm{u}_h(\bm{x}) = \sum^{n}_{i=1} B^i(\bm{x}) \bm{u}_i\,,\quad\bm{v}_h(\bm{x}) = \sum^{n}_{i=1} B^i(\bm{x}) \bm{v}_i\,,
\end{align*}
where $\bm{u}_i,\bm{v}_i\in\mathbb{R}^3$ are their control point coefficients and $B^i=\hat{B}^i\circ\bm{F}^{-1}$.
Plugging them into \eqref{eq:a} the bilinear form becomes
\begin{align}
  a(\bm{v}_h,\,\bm{u}_h) = \int_\Omega \nabla^s\bm{v}_h(\bm{x}):\mathbb{C}(\bm{x}):\nabla^s\bm{u}_h(\bm{x})\diff\bm{x}\,,
\end{align}
that can be expressed as:
\begin{align} \label{eq:a_0}
  a(\bm{v}_h,\,\bm{u}_h) = \sum_{i,j=1}^{n} \bm{v}_i \cdot \bm{K}^{ij} \bm{u}_j =\bm{\mathsf{v}}^\top \bm{\mathsf{K}} \bm{\mathsf{u}}\,,
\end{align}
where  $\bm{\mathsf{K}}$ is the problem's stiffness matrix, $\bm{\mathsf{u}}^\top = [{\bm{u}_1}^\top,{\bm{u}_2}^\top,\dots,{\bm{u}_n}^\top]^\top$ and $\bm{\mathsf{v}}$ is built in the same way as $\bm{\mathsf{u}}$.
Using Voigt's notation, the matrices $\bm{K}^{ij}\in\mathbb{R}^{3\times3}$ can be computed as (see, e.g., \cite{hughes_finite_1987,zienkiewicz_finite_2013})
\begin{align}\label{eq:K_1}
  \bm{K}^{ij} = \int_\Omega \BB^i{^\top}(\bm{x})\, \DD(\bm{x})\, \BB^j(\bm{x})\diff\bm{x}\,,
\end{align}
where the strain-displacement matrices $\BB^i\in\mathbb{R}^{6\times3}$ are calculated as
\begin{align}\label{eq:B_0}
  \BB^i{^\top}(\bm{x}) = \begin{pmatrix}
    \nabla B^i_1(\bm{x}) & 0 & 0 & \nabla B^i_2(\bm{x}) & \nabla B^i_3(\bm{x}) & 0\\
    0 & \nabla B^i_2(\bm{x}) & 0 & \nabla B^i_1(\bm{x}) & 0 & \nabla B^i_3(\bm{x})\\
    0 & 0 & \nabla B^i_3(\bm{x}) & 0 & \nabla B^i_1(\bm{x}) & \nabla B^i_3(\bm{x})
  \end{pmatrix}\,,
\end{align}
and $\nabla B^i_k(\bm{x})$, for $k=\{1,2,3\}$, is the $k$-th component of the gradient vector $\nabla B^i(\bm{x})\in\mathbb{R}^3$.
Accordingly, the material matrix $\DD(\bm{x})\in\mathbb{R}^{6\times6}$, Voigt's representation of the $\CC(\bm{x})$, is computed as:
\begin{align}\label{eq:D}
  \DD(\bm{x}) = \begin{pmatrix}
	\CC_{1111}(\bm{x}) & \CC_{1122}(\bm{x}) & \CC_{1133}(\bm{x}) & \CC_{1112}(\bm{x}) & \CC_{1113}(\bm{x}) & \CC_{1123}(\bm{x}) \\
	           & \CC_{2222}(\bm{x}) & \CC_{2233}(\bm{x}) & \CC_{2212}(\bm{x}) & \CC_{2213}(\bm{x}) & \CC_{2223}(\bm{x}) \\
	           &            & \CC_{3333}(\bm{x}) & \CC_{3312}(\bm{x}) & \CC_{3313}(\bm{x}) & \CC_{3323}(\bm{x}) \\
	           &            &            & \CC_{1212}(\bm{x}) & \CC_{1213}(\bm{x}) & \CC_{1223}(\bm{x}) \\
	           & \text{sym.}&            &            & \CC_{1313}(\bm{x}) & \CC_{1323}(\bm{x}) \\
	           &            &            &            &            & \CC_{2323}(\bm{x})
  \end{pmatrix}\,.
\end{align}

By pulling-back the gradient $\nabla B^i(\bm{x})$ to the parametric domain through $\nabla \hat B^i = \nabla B^i \circ \bm{F}$, where
\begin{align}\label{eq:gradient}
\nabla B^i = \hat{D}\bm{F}^{-\top} \hat\nabla \hat B^i\,,
\end{align}
it is possible to compute the integral \eqref{eq:K_1} in the parametric domain $\hat\Omega$ as:
\begin{align}\label{eq:K_2}
  \bm{K}^{ij} = \int_{\hat\Omega}
	\hat{\BB}^i{^\top}(\bm{\xi})\,\hat\DD(\bm{\xi})\,\hat{\BB}^j(\bm{\xi})\,
	\lvert \hat D\bm{F}(\bm{\xi})\rvert\,\diff\bm{\xi}\,,
\end{align}
where
\begin{align}\label{eq:B_1}
  \hat{\BB}^i{^\top}(\bm{\xi}) = \begin{pmatrix}
    \nabla \hat B^i_1(\bm{\xi}) & 0 & 0 & \nabla \hat B^i_2(\bm{\xi}) & \nabla \hat B^i_3(\bm{\xi}) & 0\\
    0 & \nabla \hat B^i_2(\bm{\xi}) & 0 & \nabla \hat B^i_1(\bm{\xi}) & 0 & \nabla \hat B^i_3(\bm{\xi})\\
    0 & 0 & \nabla \hat B^i_3(\bm{\xi}) & 0 & \nabla \hat B^i_1(\bm{\xi}) & \nabla \hat B^i_3(\bm{\xi})
  \end{pmatrix}\,
\end{align}
and $\hat\DD$ is the Voigt's representation of the pulled-back tensor $\hat\CC=\CC\circ\bm{F}$.

\subsection{Geometric structure of laminated composites} \label{sec:multi-layered}
\begin{figure}
\centering
\includegraphics[width=0.95\textwidth]{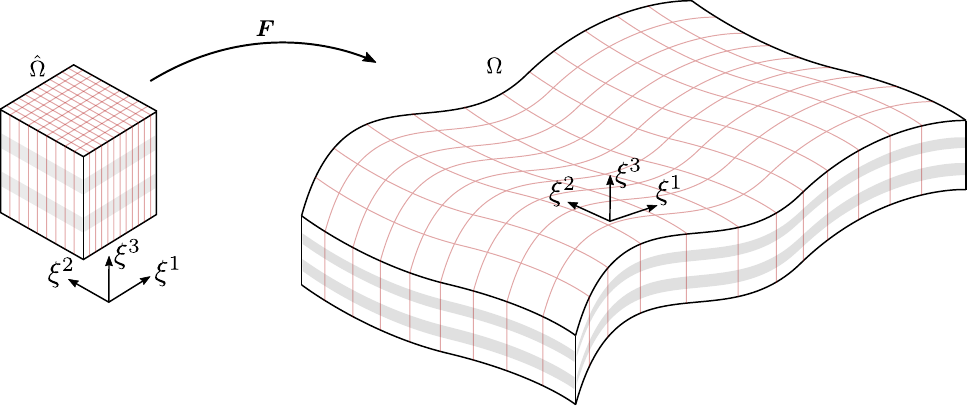}
\caption{Parametrization of a laminate composite structure. Different material layers are designated with alternated white and gray colors. Red lines refer to the underlying mesh.}
\label{fig:shell}
\end{figure}

Let us now consider the formation of stiffness matrices for 3D laminated
composite structures, that present a multi-layered structure along the shell thickness, as represented in Figure \ref{fig:shell}.
Here we assume that, as shown in that figure, the third parametric coordinate $\xi^3$ corresponds to the thickness direction.
Then, the parametric domain $\hat{\Omega}=[0,1]^3\subset\mathbb{R}^3$ can be decomposed as $\hat{\Omega} = \hat{\bar{\Omega}}\times\hat{h}$ into its in-plane $\hat{\bar{\Omega}}=[0,1]^2\subset\mathbb{R}^2$ and out-of-plane $\hat{h}=[0,1]\subset\mathbb{R}$ sub-domains.

Along the thickness direction, $m$ layers are considered, being the coordinates of the layer interfaces $\{\hat t_0, \hat t_1, \hat t_2,\dots,\hat t_{m}\}$, such that $\hat t_0 = 0$, $\hat t_1 = 1$ and $\hat t_{i+1}>\hat t_i$, for $i=0,1,\dots,m-1$.
Thus, the parametric domain of a single layer is defined as
\begin{align}
\hat{\Omega}_i = \hat{\bar{\Omega}}\times \hat{h}_i,\quad\text{for }i=1,\dots,m\,,
\end{align}
with $\hat{h}_i = [\hat{t}_{i-1},\hat{t}_i]$ and $\hat h = \hat h_1\cup\hat h_2\cup\dots\cup\hat h_m$.
By an abuse of notation we construct the image of every layer in the physical domain as $\Omega_i=\bm{F}(\hat{\Omega}_i)$.

While three parametric coordinates $\bm{\xi}=(\xi^1,\xi^2,\xi^3)$ were used for the parametric domain $\hat\Omega$,
we now split them into the in-plane coordinates $\bar{\bm{\xi}} = (\xi^1,\xi^2)$
for $\hat{\bar{\Omega}}$, and the out-of-plane coordinate $\xi^3$ for $\hat h$ (cf.\ Figure \ref{fig:shell}).

\begin{assumption} \label{assmpt:1}
We assume that the map gradient $\hat D\bm{F}$ does not depend on the third parametric direction,
but only on the in-plane ones, i.e.\ $\hat D\bm{F}=\hat D\bm{F}(\bar{\bm{\xi}})$.
\end{assumption}

This is a reasonable assumption in the case of 3D shell-like structures, whose parametrization is frequently created as the extrusion of a 2D freeform manifold along a constant direction:
\begin{align}\label{eq:extrusion}
\bm{F}(\bm{\xi}) = \bm{S}(\bar{\bm{\xi}}) + \xi^3\bm{a}\,,
\end{align}
where $\bm{S}:\mathbb{R}^2\to\mathbb{R}^3$ is the 2D manifold and $\bm{a}\in\mathbb{R}^3$ the extrusion direction.
It is easy to realize in \eqref{eq:extrusion} that the gradient $\hat D\bm{F}$ does not depend on $\xi^3$.

We consider, for every layer along the thickness, a different material.
Thus, a different elasticity tangent tensor $\hat\CC$ is associated to every layer:
\begin{align}\label{eq:C_layers}
\hat\CC(\bm{\xi}) = \left\lbrace
\begin{array}{ll}
\hat\CC_1(\bm{\xi}) & \text{if }\xi^3\in\hat h_1\,,\\
\hat\CC_2(\bm{\xi}) & \text{if }\xi^3\in\hat h_2\,,\\
\dots & \dots\\
\hat\CC_m(\bm{\xi}) & \text{if }\xi^3\in\hat h_m\,.
\end{array}
\right.
\end{align}
Then, for the case of laminated composites, the stiffness matrix \eqref{eq:K_2} can be
computed like:
\begin{align}\label{eq:Km_0}
  \bm{K}^{ij} = \sum_{l=1}^m \int_{\hat\Omega_l}
	\hat{\BB}^i{^\top}(\bm{\xi})\,\hat\DD_l(\bm{\xi})\,\hat{\BB}^j(\bm{\xi})\,
	\lvert\hat{D}\bm{F}(\bm{\xi})\rvert\,\diff\bm{\xi}\,,
\end{align}
where $\hat\DD_l$ is the Voigt's representation of every tensor $\hat\CC_l$.

\subsection{Computational cost of classical stiffness matrix assembly} \label{sec:classical_assembly}
The standard way of computing the integrals above is to calculate their contribution element by element, considering only the functions that have support in every element, and then assembling all of them together.
For that purpose, a quadrature based on the multi-layered structure must be aapplied.

Let us now, for the sake of exposition's clarity and without constituting any limitation, assume that the structure's mesh presents only one element through the thickness, as sketched in Figure \ref{fig:shell}.
This is the case of ``solid shell'' discretizations.
Thus, if the basis functions used in the discretization have degree $p$ along the three parametric directions, we consider an in-plane quadrature rule with $(p+1)^2$ points per element, while, in order to integrate precisely along the thickness, at least $p+1$ quadrature points must be used for every material layer.
Hence, the total number of quadrature points per in-plane element is $m\,(p+1)^3$.
In addition, there are $(p+1)^3$ non-vanishing basis functions in each element, therefore, the assembly's computational complexity of every in-plane element is $\mathcal{O}(m\,p^9)$ (see \cite{antolin_efficient_2015} for a further discussion).

As it can be seen, the number of quadrature points, and therefore the computational cost, scales linearly with the number of layers $m$.
When a high number of layers is considered, the assembly of the stiffness matrix becomes very expensive in terms of floating point operations.

\begin{remark}
Even if for exposition purposes we limited our discussion to the case of a single element through the thickness, the previous result can be extended in a straightforward manner to the case in which the discretization has more than one element along the third parametric direction, according to the Assumption \ref{assmpt:0}.
In this case, as long as $m\,(p+1)$ quadrature points are used along the whole laminate's thickness, independently of the discretization along the out-of-plane direction, the above complexity's estimation above is still valid.
\end{remark}

\begin{remark}
In the case of wedge elements both the number of in-plane basis functions per element and quadrature points scale as $\mathcal{O}(p^2)$.
While they present $p$ out-of-plane functions and quadrature points.
Therefore, the previous discussion regarding the method's complexity is still valid.
The same applies to the complexity results of the proposed fast method in Section \ref{sec:cost}.
\end{remark}

In the next section we propose an alternative strategy for assembling the stiffness matrix.
This new methodology presents a lower computational cost and guarantees its exactness up to machine precision.

\section{Fast assembly of stiffness matrices for 3D laminate composite structures} \label{sec:fast_assembly}
Let us first introduced the last assumption in which this work relies upon.
\begin{assumption} \label{assmpt:2}
We assume that for every single layer, the material only presents in-plane variations, and not out of plane.
I.e., the tensor $\hat\CC$ can be expressed as:
\begin{align}
\hat\CC(\bm{\xi}) = \left\lbrace
\begin{array}{ll}
\hat\CC_1(\bar{\bm{\xi}}) & \text{if }\xi^3\in\hat h_1\,,\\
\hat\CC_2(\bar{\bm{\xi}}) & \text{if }\xi^3\in\hat h_2\,,\\
\dots & \dots\\
\hat\CC_m(\bar{\bm{\xi}}) & \text{if }\xi^3\in\hat h_m\,.
\end{array}
\right.
\end{align}
\end{assumption}
In addition to purely homogeneous material plies, the assumption above is fulfilled by numerous 3D laminates, as it is the case of some fiber reinforced plies (e.g., carbon fiber) or honeycomb panels, that can present an heterogenous (and anisotropic) in-plane behavior, but their mechanical properties can be assumed homogeneous along the ply's thickness.

\newcommand{\hatS}[1]{\hat{S}^{#1_s} (\xi^1,\xi^2)}
\newcommand{\hatSp}[1]{\hat{S}^{#1_s} (\bar{\bm{\xi}})}
\newcommand{\hatT}[1]{\hat{T}^{#1_t} (\xi^3)}
\newcommand{\partialderSu}[1]{\dfrac{\partial\hat{S}^{#1_s} (\xi^1,\xi^2)}{\partial \xi^1}}
\newcommand{\partialderSv}[1]{\dfrac{\partial\hat{S}^{#1_s} (\xi^1,\xi^2)}{\partial \xi^2}}
\newcommand{\partialderT}[1]{\dfrac{\partial\hat{T}^{#1_t} (\xi^3)}{\partial \xi^3}}

We now explicitly introduce the three components of the gradient vectors $\hat\nabla\hat B^i$ present in \eqref{eq:gradient}:
\begin{align}
\hat\nabla \hat{B}^i(\bm\xi) =
\begin{pmatrix}
%
\partialderSu{i}\,\hatT{i}\\
\partialderSv{i}\,\hatT{i}\\
\hatS{i}\,\partialderT{i}
\end{pmatrix}\,.
\end{align}

Based on Assumptions \ref{assmpt:0} and \ref{assmpt:1}, the basis function gradients can be split in their in-plane and out-of-plane components as:
\begin{align}
\nabla \hat{B}^i(\bm{\xi}) = \hat{D}\bm{F}^{-\top}(\bar{\bm{\xi}}) \left(
\hat\nabla \hat{V}^{i_s}(\bar{\bm{\xi}})\,\hatT{i} + \hat{W}^{i_s}(\bar{\bm{\xi}})\,\partialderT{i}\right)\,,
\end{align}
where
\begin{align*}
\hat\nabla \hat{V}^{i_s}(\bar{\bm{\xi}})=
\begin{pmatrix}
\partialderSu{i}\\
\partialderSv{i}\\
0
\end{pmatrix}\,,\quad
 \hat{W}^{i_s}(\bar{\bm{\xi}})=
\begin{pmatrix}
0 \\ 0 \\
\hatSp{i}
\end{pmatrix}\,.
\end{align*}
Then, $\nabla \hat{B}^i$ can be expressed as
\begin{align}
\nabla \hat{B}^i(\bm{\xi}) = 
\nabla \hat{V}^{i_s}(\bar{\bm{\xi}})\,\hatT{i} + W^{i_s}(\bar{\bm{\xi}})\,\partialderT{i}\,,
\end{align}
with
\begin{subequations}
\begin{align}
\nabla \hat{V}^{i_s}(\bar{\bm{\xi}}) &= \hat{D}\bm{F}^{-\top}(\bar{\bm{\xi}})\,\hat\nabla \hat{V}^{i_s}(\bar{\bm{\xi}})\,,\\
W^{i_s}(\bar{\bm{\xi}}) &= \hat{D}\bm{F}^{-\top}(\bar{\bm{\xi}})\,\hat{W}^{i_s}(\bar{\bm{\xi}})\,.
\end{align}
\end{subequations}
Thus, the strain-displacement matrix $\hat\BB^i$ in \eqref{eq:B_1} can be also split in its in-plane and out-of-plane contributions:
\begin{align}\label{eq:B_2}
\hat\BB^i(\bm{\xi}) = \hat{\BB}^{i_s}_1(\bar{\bm{\xi}})\,\hatT{i} + \hat{\BB}^{i_s}_2(\bar{\bm{\xi}})\,\partialderT{i}\,,
\end{align}
where the in-plane matrices $\hat{\BB}^{i_s}_1$ and $\hat{\BB}^{i_s}_2$ are built by substituting the gradient $\nabla\hat B^i(\bm x)$ in \eqref{eq:B_1} with $\nabla \hat{V}^{i_s}(\bar{\bm{\xi}})$ and $W^{i_s}(\bar{\bm{\xi}})$, respectively.
Hence, considering Assumptions \ref{assmpt:1} and \ref{assmpt:2} and plugging the splitting \eqref{eq:B_2} of $\hat\BB^i$ into the computation of the stiffness matrix \eqref{eq:Km_0}, we obtain:
\begin{align}\label{eq:Km_1}
\begin{split}
  \bm{K}^{ij} = \sum_{l=1}^{m} \int_{\hat h_l} \int_{\hat{\bar\Omega}}
  \left(\hat{\BB}^{i_s}_1(\bar{\bm{\xi}})\,\hatT{i} + \hat{\BB}^{i_s}_2(\bar{\bm{\xi}})\,\partialderT{i}\right)
  \,\hat\DD_l(\bar{\bm{\xi}})\\
  \left(\hat{\BB}^{j_s}_1(\bar{\bm{\xi}})\,\hatT{j} + \hat{\BB}^{j_s}_2(\bar{\bm{\xi}})\,\partialderT{j}\right)
	\lvert\hat{D}\bm{F}(\bar{\bm{\xi}})\rvert\,\diff\bar{\bm{\xi}}\,\diff\xi^3\,.
\end{split}
\end{align}
Finally, gathering the in-plane and the out-of-plane terms, the previous expression can be reformulated as
\begin{align}\label{eq:Km_2}
\begin{split}
  \bm{K}^{ij} = \sum_{l=1}^{m}\left(
  \bm{\mathsf{P}}^{{ij_s}}_{l,11}\,\mathsf{Q}^{ij_t}_{l,11}
+ \bm{\mathsf{P}}^{{ij_s}}_{l,12}\,\mathsf{Q}^{ij_t}_{l,12}
+ \bm{\mathsf{P}}^{{ij_s}}_{l,21}\,\mathsf{Q}^{ij_t}_{l,21}
+ \bm{\mathsf{P}}^{{ij_s}}_{l,22}\,\mathsf{Q}^{ij_t}_{l,22}\right)\,,
\end{split}
\end{align}
where
\newcommand{\Poperator}{\bm{\mathsf{P}}^{{ij_s}}_{l,\alpha \beta}}
\begin{align}\label{eq:P}
\Poperator = \int_{\hat{\bar\Omega}}
\hat{\BB}^{i_s}_\alpha(\bar{\bm{\xi}})\,\hat\DD_l(\bar{\bm{\xi}})\,\hat{\BB}^{j_s}_\beta(\bar{\bm{\xi}})\,
\lvert\hat{D}\bm{F}(\bar{\bm{\xi}})\rvert\,\diff\bar{\bm{\xi}}\,,\text{ for }\alpha,\beta=\{1,2\}\,,
\end{align} 
and:
\begin{subequations}\label{eq:Q}
\begin{align}
\mathsf{Q}^{ij_t}_{l,11} &= \int_{\hat h_l} \hatT{i}\, \hatT{j} \diff\xi^3\,,\\
\mathsf{Q}^{ij_t}_{l,12} &= \int_{\hat h_l} \partialderT{i}\, \hatT{j} \diff\xi^3\,,\\
\mathsf{Q}^{ij_t}_{l,21} &= \int_{\hat h_l} \hatT{i}\, \partialderT{j} \diff\xi^3\,,\\
\mathsf{Q}^{ij_t}_{l,22} &= \int_{\hat h_l} \partialderT{i}\, \partialderT{j} \diff\xi^3\,.
\end{align}
\end{subequations}

Thus, the 3D integrals involved in the computation of the stiffness matrix \eqref{eq:Km_0} are decomposed in \eqref{eq:Km_2} as combinations of 2D \eqref{eq:P} and 1D \eqref{eq:Q} integrals.
It is worth noting that the stiffness matrices computed by applying \eqref{eq:Km_2} will be identical, up to machine precision, to the ones obtained with \eqref{eq:Km_0}, no approximation is introduced.

The operators $\Poperator$ are computed as 2D integrals for every material layer, however, from layer to layer, only the material
term $\hat\DD_l$ changes. Thus, in the cases in which the same material is used with different orientations, the 
matrices $\Poperator$ are needed to be computed only once for each of them. E.g., in the classical cross-ply design in which the same fibered material is stacked with $0^\circ/90^\circ/0^\circ/90^\circ/\dots$ orientations, the operators
$\Poperator$ will be computed only for two different layers, one corresponding to the $0^\circ$ orientation, and another one for $90^\circ$.
In the also common case of $0^\circ/\pm45^\circ/90^\circ/\dots$ orientations, the matrices $\Poperator$ will be computed for four different cases only.

\begin{remark} \label{rmk:orientations}
As it is detailed in, e.g.,  \cite[Section 3.3-4]{reddy_practical_1995}, the material matrix $\DD$ of an orthotropic material (satisfying the Assumption \ref{assmpt:2}) with a given in-plane orientation angle $\theta$ can be split
as:
\begin{align}
\DD = \bm{\mathsf{A}}_1 + \cos^4\theta\bm{\mathsf{A}}_2 + \cos^3\theta\sin\theta\bm{\mathsf{A}}_3 + \cos^2\theta\bm{\mathsf{A}}_4 + \cos\theta\sin\theta\bm{\mathsf{A}}_5\,,
\end{align}
where the matrices $\bm{\mathsf{A}}_i\in\mathbb{R}^{6\times6}$ only depend on the material properties, but not on the orientation $\theta$.
The same decomposition applies to $\CC$.

Therefore, in the frequent case in which the same material is used for the different layers, but changing its orientation from layer to layer, a maximum of five different matrices $\Poperator$, one for every operator $\bm{\mathsf{A}}_i$, must be computed.

On the other hand, in \ref{ap:alternative} we propose an alternative way for computing the $\Poperator$ matrices without the use of Voigt's notation.
\end{remark}


\subsection{Computational cost of fast stiffness matrix assembly} \label{sec:cost}
Studying Equations \eqref{eq:P} and \eqref{eq:Q} it is simple to realize that the computational cost of $\mathsf{Q}^{ij_t}_{l,\alpha\beta}$ is negligible compared
to the cost of $\Poperator$.
In fact, $\mathsf{Q}^{ij_t}_{l,\alpha\beta}$ can be even pre-computed analytically, without the use of a numerical quadrature.
Thus, the complexity of the proposed assembly technique is bounded by the computational cost of the matrices $\Poperator$.

The computation of $\Poperator$ is carried out by assembling its contribution for every in-plane element.
As before, considering degree $p$ along all the parametric directions, $(p+1)^2$ non-zero in-plane basis functions and $(p+1)^2$ in-plane quadrature points must be used in the integration of every single in-plane element.
Therefore, the computational complexity of every matrix $\Poperator$, considering all the involved layers, is $\mathcal{O}(\bar{m}\,p^6)$,
where $\bar{m}$ is the number of different $\DD_l$ operators considered in the laminated structure.
E.g., $\bar m=2$ for the stack design $0^\circ/90^\circ/0^\circ/90^\circ/\dots$, and
$\bar m = 4$ for $0^\circ/\pm45^\circ/90^\circ/\dots$.
As stated in the Remark \ref{rmk:orientations}, for the common case in which the same anisotropic material is used, but stacked with different orientations, it holds $\bar m \leq 5$.

Thus, the computational complexity of the proposed method is much smaller than the one of the standard integration procedure: $\mathcal{O}(\bar{m}\,p^6)$ against $\mathcal{O}(m\,p^9)$, with $\bar m \leq m$ (cf.\ Section \ref{sec:classical_assembly}).
Additionally, for a small value $\bar m$ and a high number of layers, it holds $\bar m << m$, what makes the difference between both complexities even higher.

\begin{remark}
As far as the computational cost of the stiffness matrix assembly is bounded by the calculation of the 2D operators $\Poperator$, the assembly complexity is independent of the number of elements or functions used along the thickness.
Therefore, it is possible to use any degree, number of elements or continuity for the discretization along thickness at the same computational cost for the assembly.

This can be useful, for instance, in the case of solid shells discretizations, where the use of more than one element through the thickness would alleviate possible locking phenomena.

Moreover, our approach makes the use of 3D models viable: as pointed out in \cite{liew_overview_2019}, the assembly burden is one the main obstacles that prevents their use in the analysis of laminated composites with a large number of laminae.
However, it is worth reminding that other operations, such as the solution of the linear system of equations, are still dominated by the number of degrees of freedom and function's continuity (see, e.g., the discussion in \cite{collier_cost_2012,collier_cost_2013} for isogeometric discretizations).
\end{remark}

\begin{remark}
In the computation of $\Poperator$, in Equation \eqref{eq:P}, it is possible to apply the sum-factorization technique proposed in \cite{antolin_efficient_2015}, that would allow to reduce the computation complexity to $\mathcal{O}(\bar{m}\,p^5)$.
This improvement has not been detailed in this work for the sake of brevity.
By directly applying the sum-factorization method as described in \cite{antolin_efficient_2015}, without splitting the assembly in its in-plane/out-of-plane contributions, a theoretical computational complexity $\mathcal{O}(m\,p^7)$ is expected.

In addition, considering the weighted-quadrature techniques proposed in \cite{calabro_fast_2017,hiemstra_fast_2019} it would be potentially possible to improve, even further, the dependency on $p$ of the computational complexity of $\Poperator$ in the case of highly continuous spline basis functions.
Another possibility for reducing that cost would be to apply the low-rank approximations techniques, as proposed in \cite{mantzaflaris_low_2017} also in the context of isogeometric analysis.
\end{remark}

\section{Numerical experiment} \label{sec:numerical}
In this section we aim at illustrating the performance of the proposed assembly method (compared to the standard one) by computing the stiffness matrix associated of the classical Pagano plate problem \cite{pagano_exact_1970} in the context of tensor-product isogeometric discretizations.
This problem consists in a flat rectangular plate composed of a group of stacked layers of the same linear elastic orthotropic material distributed with different orientations.
The material mechanical properties are: $E_1 = 25\,\text{GPa}$, $E_2=E_3 = 1\,\text{GPa}$, $\mu_{12} = \mu_{13} = 0.2\,\text{GPa}$, $\mu_{23} = 0.5\,\text{GPa}$ and $\nu_{12} = \nu_{13} = \nu_{23} = 0.25$, and the material layers are stacked in two different configurations, $0^\circ/90^\circ/0^\circ/90^\circ/\dots$ and $0^\circ/\pm45^\circ/90^\circ/\dots$.

The same discretization degree $p$ is used along all the parametric directions, and different values of $p$ are considered.
Regarding the mesh discretization, different number of elements along the in-plane directions are considered, while one element through the thickness is used.
In all the cases, non-rational $C^{p-1}$ continuous B-splines are used.

All the results shown below were obtained with implementations of the standard assembly method and the fast one proposed in this work based on the isogeometric analysis library GeoPDEs version 3.0 \cite{vazquez_new_2016}.

\begin{figure}
 \center
 \subfigure[$p=1$\label{fig:res1_1}]{\ifrecompiletikz\else\includegraphics{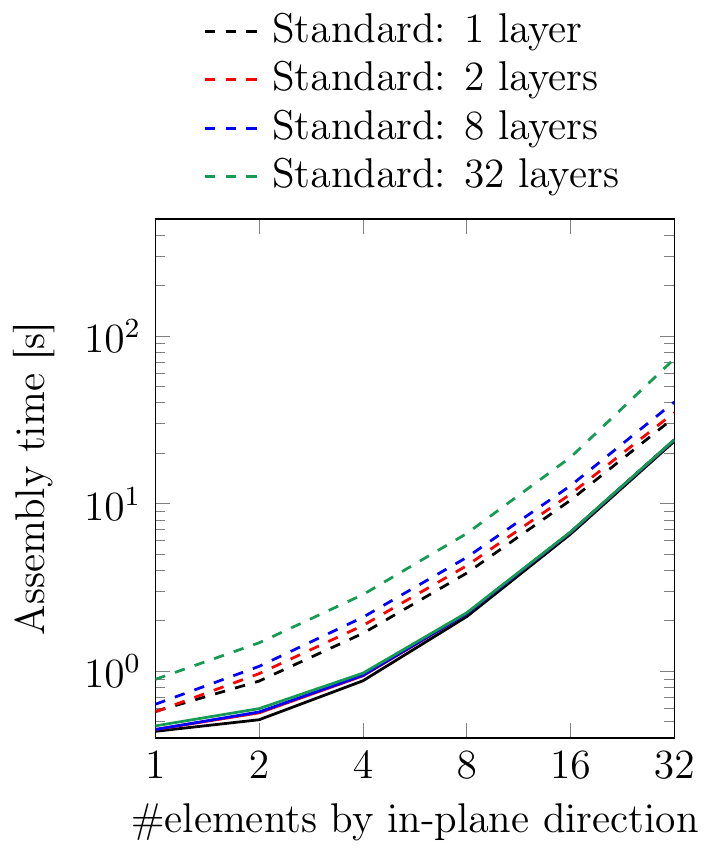}\fi}
 \subfigure[$p=2$\label{fig:res1_2}]{\ifrecompiletikz\else\includegraphics{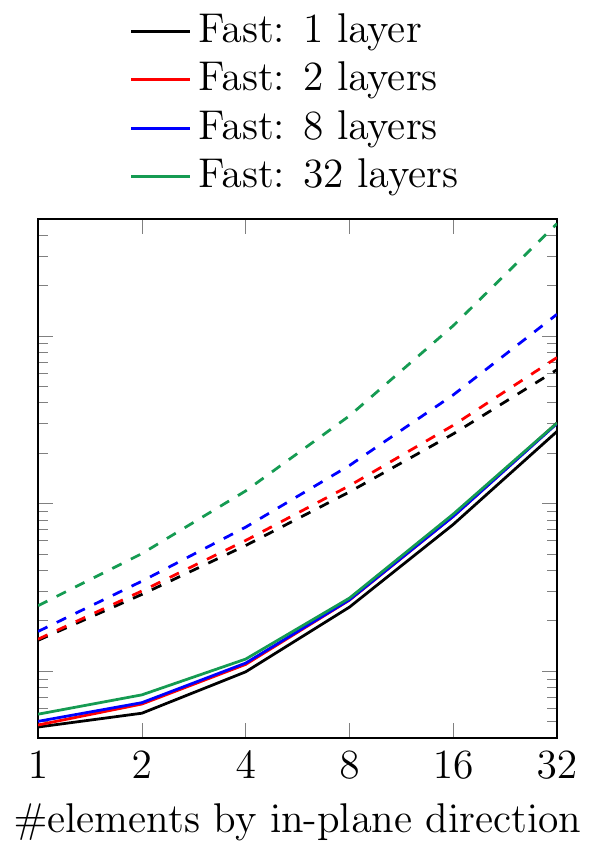}\fi}\\
 \subfigure[$p=3$\label{fig:res1_3}]{\ifrecompiletikz\else\includegraphics{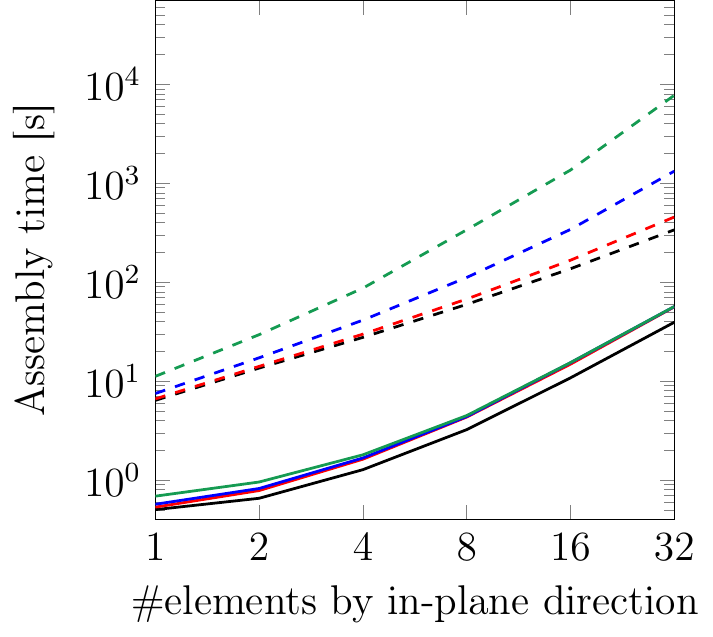}\fi}
 \subfigure[$p=4$\label{fig:res1_4}]{\ifrecompiletikz\else\includegraphics{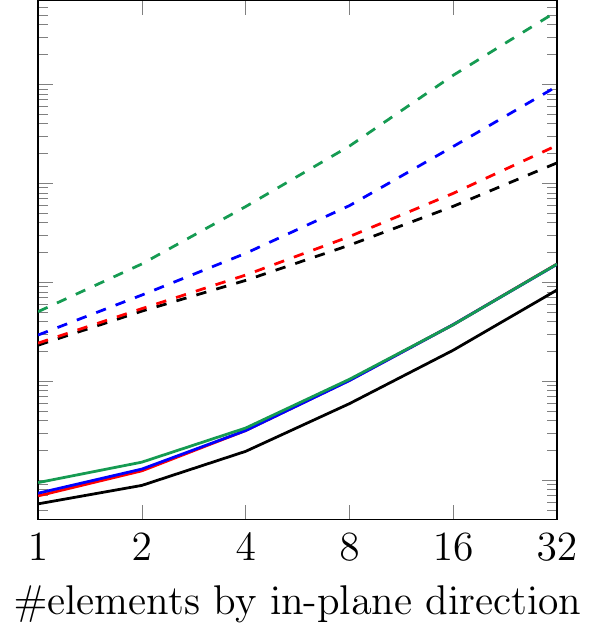}\fi}\\
 \caption{Pagano plate problem: $0^\circ/90^\circ/0^\circ/90^\circ\dots$ cross-ply laminated composite.
 Comparison of assembly times between the standard  procedure and the fast one proposed in this work,
 respect to the number of elements by in-plane direction for different degrees and numbers of material layers.}
 \label{fig:res1}
\end{figure}

\begin{figure}
 \center
 \subfigure[$p=1$\label{fig:res2_1}]{\ifrecompiletikz\else\includegraphics{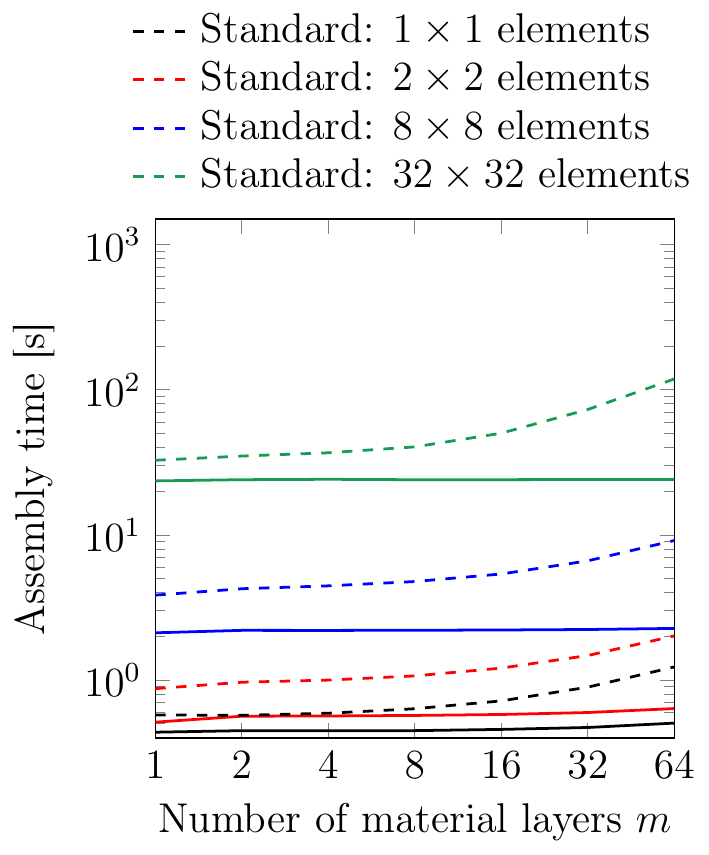}\fi}
 \subfigure[$p=2$\label{fig:res2_2}]{\ifrecompiletikz\else\includegraphics{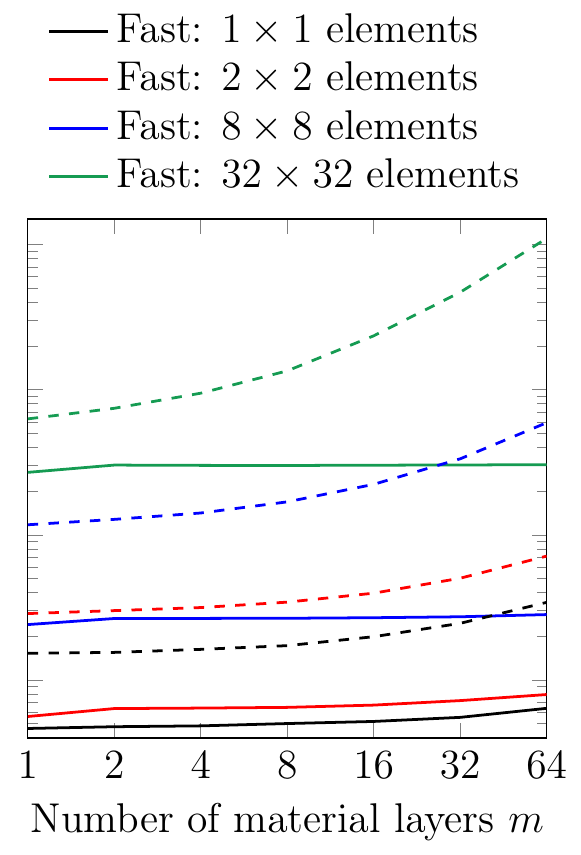}\fi}\\
 \subfigure[$p=3$\label{fig:res2_3}]{\ifrecompiletikz\else\includegraphics{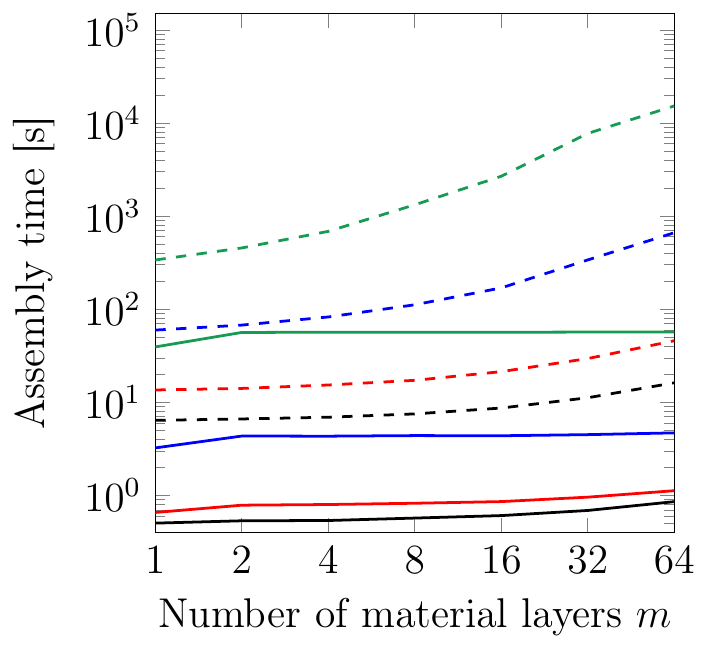}\fi}
 \subfigure[$p=4$\label{fig:res2_4}]{\ifrecompiletikz\else\includegraphics{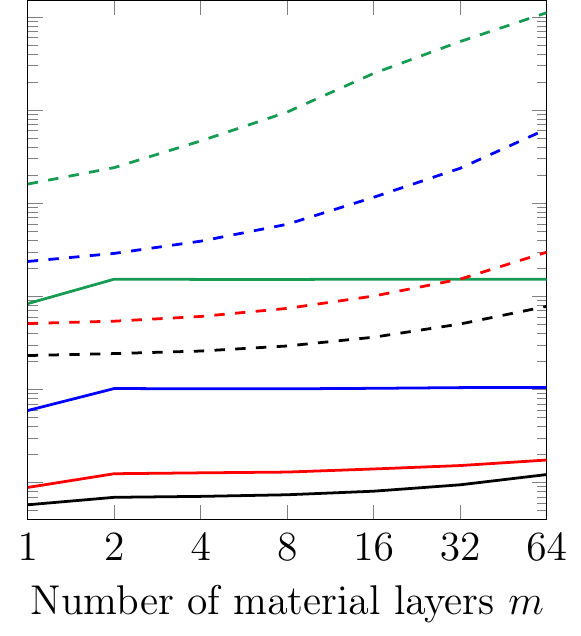}\fi}\\
 \caption{Pagano plate problem: $0^\circ/90^\circ/0^\circ/90^\circ\dots$ cross-ply laminated composite.
 Comparison of assembly times between the standard  procedure and the fast one proposed in this work,
 respect to the number of material layers for different discretizations.}
 \label{fig:res2}
\end{figure}

\begin{figure}
 \center
 \subfigure[$p=1$\label{fig:res3_1}]{\ifrecompiletikz\else\includegraphics{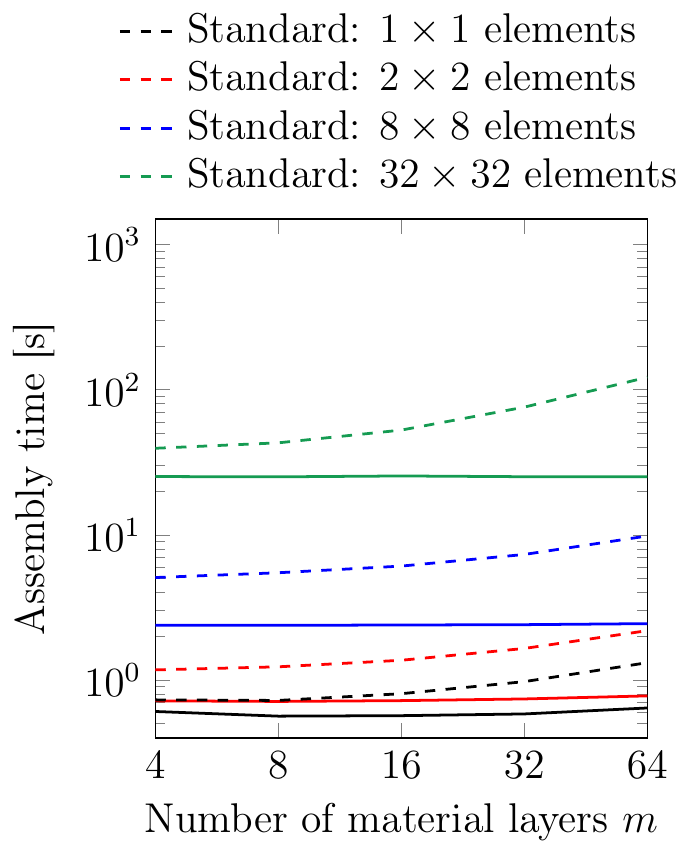}\fi}
 \subfigure[$p=2$\label{fig:res3_2}]{\ifrecompiletikz\else\includegraphics{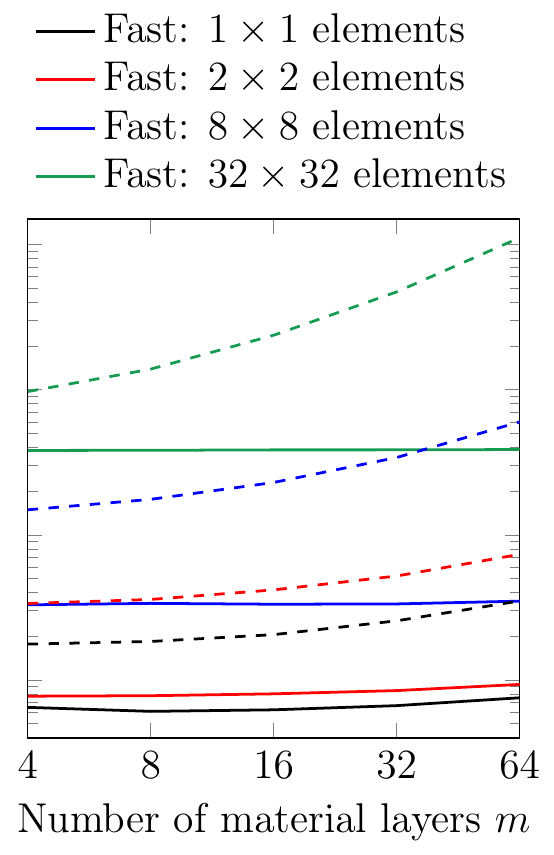}\fi}\\
 \subfigure[$p=3$\label{fig:res3_3}]{\ifrecompiletikz\else\includegraphics{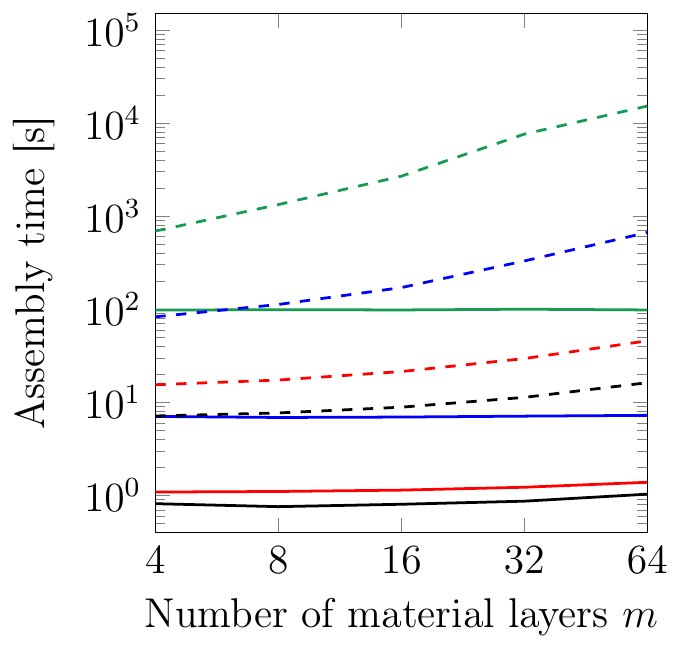}\fi}
 \subfigure[$p=4$\label{fig:res3_4}]{\ifrecompiletikz\else\includegraphics{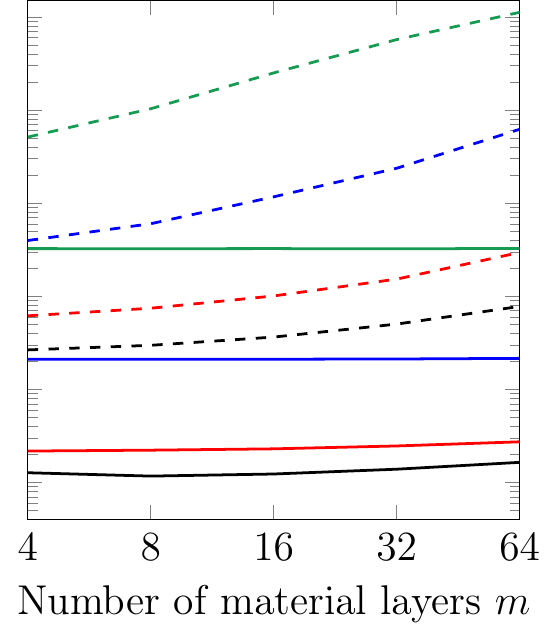}\fi}\\
 \caption{Pagano plate problem: $0^\circ/\pm45^\circ/90^\circ/\dots$ cross-ply laminated composite.
 Comparison of assembly times between the standard procedure and the fast one proposed in this work,
 respect to the number of material layers for different discretizations.}
 \label{fig:res3}
\end{figure}

In Figures \ref{fig:res1} and \ref{fig:res2} the assembly times of the standard method
and the fast one proposed in this work are compared for the $0^\circ/90^\circ/0^\circ/90^\circ/\dots$ configuration.
Figure \ref{fig:res1} shows the assembly times as a function of the number of elements along the in-plane direction, for different numbers of materials
and degrees; and, on the other hand, the assembly times as a function of the number of layers are shown in Figure \ref{fig:res2}.
Analogously, the dependence of the assembly times respect to the number of layers is shown in Figure \ref{fig:res3} for the $0^\circ/\pm45^\circ/90^\circ/\dots$ configuration and different discretizations.

As it can be seen in Figures \ref{fig:res1} to \ref{fig:res3}, the fast assembly method (solid lines)
outperforms the standard one (dashed lines) for all the degrees, meshes and number of layers considered.
For high degree and large number of layers, the proposed method reduces the assembly time of the standard procedure by more than two orders of magnitude.
It is also worth mentioning that the matrices computed using the standard and fast methods are equal, term by term, up to machine precision.

As shown in Figures \ref{fig:res2} and \ref{fig:res3}, and explained in Section \ref{sec:cost}, the assembly time of the standard method scales linearly with the number of layers $m$, whereas for the fast method is constant and independent of the number of layers.
The latter being true as far as the number of different material configurations $\bar m$ is $\bar m \leq m$.
E.g., as it can be appreciated in Figure \ref{fig:res2}, for the $0^\circ/90^\circ/0^\circ/90^\circ\dots$ laminate design ($\bar m = 2$), the $1$ layer case ($m = 1$) is faster than the $2$ layers case ($m = 2$): the assembly time initially scales with the number of layers.
Nevertheless, this is no longer true for $m>2$, as the condition $\bar m < m$ holds.

In addition, it is worth mentioning that in Figure \ref{fig:res2} the assembly time for $m=2$ is less than twice the time for $m=1$.
This effect is caused by the overhead of evaluating basis functions, allocating data structures, etc.
For a single layer this time overhead is not negligible respect to the computing time for evaluating the matrices $\Poperator$. This is even clearer for $p=1$ and $p=2$, where the cost of $\Poperator$ is negligible respect to other costs.

On the other hand, as it can be observed in Figures \ref{fig:res2} and \ref{fig:res3}, for small number of elements ($1\times1$ and  $2\times2$ meshes, black and red solid lines, respectively)
the assembly time is not exactly constant: it is initially flat, but presents a slight increment as the number of layers grows.
This is due to the fact that the constant cost $\mathcal{O}(\bar{m}\,p^6)$
of the computation of the $\Poperator$ matrices (for a small number of elements) dominates for a small number of layers, whereas, as the number of layers increases,
the cost of combining those matrices for $m$ layers (Equation \eqref{eq:Km_2}), with $m>>\bar m$, starts to be dominant.


Finally, we would like to remark that, as highlighted in Section \ref{sec:intro}, the proposed method, considering multiple laminae, is even faster than a simplified 3D model that uses a single layer with homogeneized material properties, assembled with the standard procedure.
The assembly times of the simplified model are equivalent to the ones obtained for the one layer cases shown in the numerical experiments (i.e., dashed black lines in Figure \ref{fig:res1} and first point of all dashed lines in Figure \ref{fig:res2}).

%

\section{Conclusions} \label{sec:conclusions}

In this paper we present a new method for assembling stiffness matrices of 3D laminate composite structures, in a linear elasticity framework, using isogeometric and finite element discretizations.
This method relies on three main assumptions: the basis functions can be multiplicatively split into its in-plane and out-of-plane components; the gradient of the geometric parametrization map is constant along the thickness direction (as it is the case of geometries built as extrusion of 2D freeforms); and the mechanical properties of each material layer can only vary in-plane, and not out-of-plane (this is the case of many materials, e.g., carbon fiber plies, that present a 2D in-plane reinforcement, being homogeneous along the ply thickness).
Based on these assumptions, we develop an assembly method in which the 3D stiffness matrix is computed as a combination of 2D and 1D integrals and that results in a matrix identical, up to machine precision, to the one obtained using the standard assembly procedure.
The theoretical assembly cost of the proposed method per in-plane element is $\mathcal{O}(\bar m\,p^6)$, where $p$ is the discretization degree and $\bar m$ the number of different material configuratioins used in the cross-ply laminated composite.
$\bar m$ is such that $\bar m \leq m$, and in general it holds $\bar m << m$, where $m$ is the number of layers in the structure. Thus, the cost of the proposed method is virtually independent of $m$.
On the other hand, the standard integration method, in which $p+1$ quadrature points along the thickness for every material layer are considered, presents a computational cost per in-plane element $\mathcal{O}(m\,p^9)$, that scales linearly with $m$.

The performance of the proposed method is illustrated by means of the classical Pagano plate problem in the context of isogeometric discretizations, in which we confront the obtained assembly times against the ones obtained with the standard method.
The results show how the proposed method outperforms the standard one for any considered discretization: for high values of $p$ and $m$ a difference higher than two orders of magnitude, in terms of assembly times between both methods, was observed.
A further improvement of the proposed method's cost, in the context of isogeometric discretizations, may be achieved by applying weighted-quadrature techniques \cite{calabro_fast_2017,hiemstra_fast_2019} for the computation of the involved 2D integrals.

\section*{Acknowledgements}
The author gratefully acknowledges the support of the European Research Council, through the ERC AdG n.\ 694515 - CHANGE grant, and also thanks Annalisa Buffa for her helpful insights and suggestions.

\appendix

\section{Fast assembly of stiffness matrices for laminated composites without using Voigt's notation} \label{ap:alternative}
Planas et al.\ introduced in \cite{planas_b_2012} an alternative methodology for
the assembly of stiffness matrices in the context of solid mechanics, denoted as ``$\BB$ free'', that overcomes the use of Voigt's notation.
Based on that approach, we detail in this Appendix the construction of the operators
$\Poperator$ defined in \eqref{eq:P} avoiding the use of Voigt's notation, i.e, without building the material matrix $\DD$ \eqref{eq:D} and the strain-displacement operators $\hat\BB^i$ \eqref{eq:B_2}.

Using the bracket operator $\bullet\left\lbrace\bullet,\bullet\right\rbrace:\mathbb{
R}^{3\times3\times3\times3}\times\mathbb{R}^3\times\mathbb{R}^3\to\mathbb{R}^{3\times3}$ defined in \cite{planas_b_2012}, the stiffness matrix \eqref{eq:K_1} can be computed as
\begin{align}
  \bm{K}^{ij} = \int_\Omega \CC(\bm{x})\left\lbrace\nabla B^i(\bm{x}), \nabla B^j(\bm{x})\right\rbrace\diff\bm{x}\,,
\end{align}
where the contraction $\CC(\bm{x})\left\lbrace\nabla B^i(\bm{x}), \nabla B^j(\bm{x})\right\rbrace\in\mathbb{R}^{3\times3}$ is such that the vector
$\nabla B^i$ is contracted with the second component of $\CC$, and $\nabla B^j$ with the fourth one, i.e.:
\begin{align}
\bm{a}\cdot\left(\CC\{\bm{b},\bm{d}\}\,\bm{c}\right) = \left(\bm{a}\otimes\bm{b}\right):\CC:\left(\bm{c}\otimes\bm{d}\right)\,.
\end{align}
Thus, being $\{\bm{e}_1,\bm{e}_2,\bm{e}_3\}$ the Cartesian orthonormal basis, the contraction
$\CC\left\lbrace\bm{e}_i,\bm{e}_j\right\rbrace\in\mathbb{R}^{3\times3}$ for small strain linear isotropic materials can be computed as \cite{planas_b_2012}:
\begin{align}
\CC\left\lbrace\bm{e}_i,\bm{e}_j\right\rbrace
= \lambda\,\bm{e}_i\otimes\bm{e}_j + \mu\left(\bm{e}_i\cdot\bm{e}_j\,\one + \bm{e}_j\otimes\bm{e}_i\right)\,,
\end{align}
where $\lambda$ and $\mu$ are the Lam\'e coefficients and $\bm{I}\in\mathbb{R}^{3\times 3}$ is the identity tensor.
For orthotropic materials the contraction reads:
\begin{align}
\begin{split}
\CC&\left\lbrace\bm{e}_i,\bm{e}_j\right\rbrace
= \lambda\,\bm{e}_i\otimes\bm{e}_j + \mu\left(\bm{e}_i\cdot\bm{e}_j\,\one + \bm{e}_j\otimes\bm{e}_i\right)\\
&+\alpha_1\bar{\bm{a}}_{1,i}\otimes\bar{\bm{a}}_{1,j}+\alpha_2\bar{\bm{a}}_{2,i}\otimes\bar{\bm{a}}_{2,j}\\
&+\alpha_3\left(\right(\bm{e}_i\cdot\bar{\bm{a}}_{1,j} + \bm{e}_i\cdot\bm{e}_j\left)\one +\bm{e}_j\otimes\bar{\bm{a}}_{1,i} + \bar{\bm{a}}_{1,j}\otimes\bm{e}_i\right)\\
&+\alpha_4\left(\right(\bm{e}_i\cdot\bar{\bm{a}}_{2,j} + \bm{e}_i\cdot\bm{e}_j\left)\one +\bm{e}_j\otimes\bar{\bm{a}}_{2,i} + \bar{\bm{a}}_{2,j}\otimes\bm{e}_i\right)\\
&+\alpha_5\left(\bm{e}_i\otimes\bar{\bm{a}}_{1,j} + \bar{\bm{a}}_{1,i}\otimes\bm{e}_j\right)+\alpha_6\left(\bm{e}_i\otimes\bar{\bm{a}}_{2,j} + \bar{\bm{a}}_{2,i}\otimes\bm{e}_j\right)\\
&+\alpha_7\left(\bar{\bm{a}}_{1,i}\otimes\bar{\bm{a}}_{2,j}+\bar{\bm{a}}_{2,i}\otimes\bar{\bm{a}}_{1,j}\right)\,,
\end{split}
\end{align}
with $\bar{\bm{a}}_{\beta,i} = \left(\bm{a}_\beta\cdot\bm{e}_i\right)\,\bm{a}_\beta$, where $\bm{a}_1$ and $\bm{a}_2$ are the main orthonormal in-plane material directions,
and $\lambda$, $\mu$, $\alpha_k$, with $k=1,\dots,7$, are the nine material coefficients (see, e.g., \cite[Section 3.3]{schroder_simple_2002} for further details).

Let us now write the map gradient $\hat D \bm{F}$ as a function of the covariant basis $\{\bm{g}_1,\bm{g}_2,\bm{g}_3\}$:
\begin{align}
\hat D \bm{F}(\bar{\bm{\xi}}) = \sum_{i=1}^{3}\bm{e}_i\otimes\bm{g}_i(\bar{\bm{\xi}})\,,
\end{align}
where the Assumption \ref{assmpt:1} was considered.
In the same way, its inverse $\hat D \bm{F}^{-\top}$ can be expressed as
\begin{align}
\hat D \bm{F}^{-\top}(\bar{\bm{\xi}}) = \sum_{i=1}^{3}\bm{e}_i\otimes\bm{g}^i(\bar{\bm{\xi}})\,,
\end{align}
where $\{\bm{g}^1,\bm{g}^2,\bm{g}^3\}$ is the contravariant basis, such that $\bm{g}^i\cdot\bm{g}_i=1$
and $\bm{g}^i\cdot\bm{g}_j=0$ if $i\neq j$, for $i,j=\lbrace 1,2,3\rbrace$.

Pulling-back the computation of $\bm{K}^{ij}$ to the parametric domain, as in \eqref{eq:K_2}, we obtain:
\begin{align}
  \bm{K}^{ij} = \int_{\hat\Omega}
  \hat\CC(\bm{\bm{\xi}})\left\lbrace\hat D\bm{F}^{-\top}(\bar{\bm{\xi}})\hat \nabla \hat B^i(\bm{\xi}),\, \hat D\bm{F}^{-\top}(\bar{\bm{\xi}})\hat \nabla \hat B^j(\bm{\xi})\right\rbrace
	\lvert\hat{D}\bm{F}(\bar{\bm{\xi}})\rvert\,\diff\bm{\xi}\,,
\end{align}
that can be rewritten as
\begin{align}\label{eq:K_3}
  \bm{K}^{ij} = \int_{\hat\Omega}
  \tilde\CC(\bm{\bm{\xi}})\left\lbrace\hat \nabla \hat B^i(\bm{\xi}), \hat \nabla \hat B^j(\bm{\xi})\right\rbrace
	\lvert\hat{D}\bm{F}(\bar{\bm{\xi}})\rvert\,\diff\bm{\xi}\,,
\end{align}
where $\tilde\CC$ is the pull-back of $\hat\CC$ with $\hat D\bm{F}^{-\top}$ for the second and fourth components.

On the other hand, $\hat\nabla\hat B^i$ can be split according to its Cartesian components as (see Assumption \ref{assmpt:0}):
\begin{align}
\begin{split}
\hat\nabla\hat B^i(\xi^1,\xi^2,\xi^3) =& \partialderSu{i}\,\hatT{i}\,\bm{e}_1 +\partialderSv{i}\,\hatT{i}\,\bm{e}_2\\
&+\hatS{i}\,\partialderT{i}\,\bm{e}_3\,,
\end{split}
\end{align}
and the contraction $\tilde\CC(\bm{\bm{\xi}})\left\lbrace\hat \nabla \hat B^i(\bm{\xi}), \hat \nabla \hat B^j(\bm{\xi})\right\rbrace\in\mathbb{R}^{3\times 3}$ becomes:
\begin{align} \label{eq:contraction_2}
\begin{split}
\tilde\CC(\bm{\bm{\xi}})&\left\lbrace\hat \nabla \hat B^i(\bm{\xi}), \hat \nabla \hat B^j(\bm{\xi})\right\rbrace =\\
&\left[
\tilde{\bm{C}}_{11}(\bm{\xi}) \partialderSu{i}\,\partialderSu{j}\right.\\
&+\tilde{\bm{C}}_{12}(\bm{\xi})\partialderSu{i}\,\partialderSv{j}\\
&+\tilde{\bm{C}}_{21}(\bm{\xi})\partialderSv{i}\,\partialderSu{j}\\
&\left.+\tilde{\bm{C}}_{22}(\bm{\xi})\partialderSv{i}\,\partialderSv{j}\right] \hatT{i}\hatT{j}\\
&+ \left[\tilde{\bm{C}}_{13}(\bm{\xi}) \partialderSu{i}+\tilde{\bm{C}}_{23}(\bm{\xi}) \partialderSv{i}\right]\\
& \hatS{j} \hatT{i}\partialderT{j}\\
&+ \left[\tilde{\bm{C}}_{31}(\bm{\xi}) \partialderSu{j}+\tilde{\bm{C}}_{32}(\bm{\xi}) \partialderSv{j}\right]\\
&\hatS{i}  \partialderT{i}\hatT{j}\\
&+\tilde{\bm{C}}_{33}(\bm{\xi}) \hatS{i}\hatS{j} \partialderT{i}\partialderT{j}
\end{split}
\end{align}
where $\tilde{\bm{C}}_{\alpha\beta}\in\mathbb{R}^{3\times3}$ is
\begin{align} \label{eq:C_contracted}
\tilde{\bm{C}}_{\alpha\beta}(\bm{\xi}) = \sum_{i,j=1}^{3} \hat\CC(\bm{\xi})\left\lbrace\bm{e}_i, \bm{e}_j\right\rbrace\, \left(\bm{g}^i(\bar{\bm{\xi}})\cdot\bm{e}_\alpha\right)\, \left(\bm{g}^j(\bar{\bm{\xi}})\cdot\bm{e}_\beta\right)\,.
\end{align}

Finally, substituting \eqref{eq:contraction_2} into \eqref{eq:K_3}, and including the Assumptions \ref{assmpt:1} and \ref{assmpt:2}, the stiffness matrix terms can be rearranged in the same way as we did for \eqref{eq:Km_2}, and the $\Poperator$ matrices \eqref{eq:P} become
\begin{subequations}
\begin{align}
\begin{split}
\bm{\mathsf{P}}^{ij_s}_{l,11} = \int_{\hat{\bar\Omega}}
&\left[
\tilde{\bm{C}}_{11}(\bar{\bm{\xi}}) \partialderSu{i}\,\partialderSu{j}\right.\\
&+\tilde{\bm{C}}_{12}(\bar{\bm{\xi}})\partialderSu{i}\,\partialderSv{j}\\
&+\tilde{\bm{C}}_{21}(\bar{\bm{\xi}})\partialderSv{i}\,\partialderSu{j}\\
&\left.
+\tilde{\bm{C}}_{22}(\bar{\bm{\xi}})\partialderSv{i}\,\partialderSv{j}
\right]\lvert\hat{D}\bm{F}(\bar{\bm{\xi}})\rvert\,\diff\bar{\bm{\xi}}\,,
\end{split}\\
\begin{split}
\bm{\mathsf{P}}^{ij_s}_{l,12} = \int_{\hat{\bar\Omega}}
 &\left[
 \tilde{\bm{C}}_{13}(\bar{\bm{\xi}}) \partialderSu{i}+\tilde{\bm{C}}_{23}(\bar{\bm{\xi}}) \partialderSv{i}\right]\\
 &\hatS{j}\,\lvert\hat{D}\bm{F}(\bar{\bm{\xi}})\rvert\,\diff\bar{\bm{\xi}}\,,
\end{split}\\
\begin{split}
\bm{\mathsf{P}}^{ij_s}_{l,21} = \int_{\hat{\bar\Omega}}
 &\left[
 \tilde{\bm{C}}_{31}(\bar{\bm{\xi}}) \partialderSu{j}+\tilde{\bm{C}}_{32}(\bar{\bm{\xi}}) \partialderSv{j}\right]\\
 &\hatS{i}\,\lvert\hat{D}\bm{F}(\bar{\bm{\xi}})\rvert\,\diff\bar{\bm{\xi}}\,,
\end{split}\\
\begin{split}
\bm{\mathsf{P}}^{ij_s}_{l,22} = \int_{\hat{\bar\Omega}}
&\tilde{\bm{C}}_{33,l}(\bar{\bm{\xi}})\, \hat{S}^{i_s}(\bar{\bm{\xi}})\, \hat{S}^{j_s}(\bar{\bm{\xi}})
\lvert\hat{D}\bm{F}(\bar{\bm{\xi}})\rvert\,\diff\bar{\bm{\xi}}\,.
\end{split}
\end{align}
\end{subequations}
The terms $\tilde{\bm{C}}_{\alpha\beta,l}$ are
computed following \eqref{eq:C_contracted} and particularizing $\hat\CC$ for every material layer $l$, with $l=1,\dots,m$.
Therefore, the matrices $\Poperator$ can be computed, avoiding the use of Voigt's notation, by means of in-plane terms only.

\bibliographystyle{abbrv}

\bibliography{biblio}
	
\end{document}

%% file: macros.tex
\newcommand{\diff}{\text{d}}

\newcommand{\strain}{\bm{\varepsilon}}
\newcommand{\stress}{\bm{\sigma}}
\newcommand{\CC}{\mathbb{C}}
\newcommand{\BB}{\bm{\mathsf{B}}}
\newcommand{\DD}{\bm{\mathsf{D}}}
\newcommand{\one}{\bm{I}}

\def\be{\begin{equation}}
\def\ee{\end{equation}}
\def\ba{\begin{array}}
\def\ea{\end{array}}

\definecolor{dgreen}{RGB}{0,139,0}